%% file: article.tex
\documentclass[11 pt]{amsart}
\input{preamble.tex}
\input{square.tex}

\author[Ara]{Dimitri Ara}
\address{Aix~Marseille~Univ,~CNRS,~Marseille,~I2M,~Marseille,~France}
\email{dimitri.ara@univ-amu.fr}

\author[Gagna]{Andrea Gagna}
\address{Institute of Mathematics, Czech Academy of Sciences\\ \v{Z}itn\'a 25 \\115 67   Praha 1\\ Czech Republic}
\email{gagna@math.cas.cz}

\author[Ozornova]{Viktoriya Ozornova}
\address{Max Planck Institute for Mathematics, Bonn, Germany}
\email{viktoriya.ozornova@mpim-bonn.mpg.de}

\author[Rovelli]{Martina Rovelli}
\address{Department of Mathematics and Statistics, 
University of Massachusetts Amherst, MA 01003-9305
Amherst, USA
}
\email{mrovelli@umass.edu} 

\title[A categorical characterization of strong Steiner $\omega$-categories]{A categorical characterization of\\strong Steiner $\omega$-categories}


\begin{document}

\maketitle

\begin{abstract}
Strong Steiner $\omega$-categories are a class of $\omega$-categories that admit algebraic models in the form of chain complexes, whose formalism allows for several explicit computations. The conditions defining strong Steiner $\omega$-categories are traditionally expressed in terms of the associated chain complex, making them somewhat disconnected from the $\omega$-categorical intuition. The purpose of this paper is to characterize this class as the class of $\omega$-categories generated by polygraphs that satisfy a loop-freeness condition that does not make explicit use of the associated chain complex and instead relies on the categorical features of $\omega$-categories.
\end{abstract}

\tableofcontents

\section*{Introduction}

Infinite-dimensional categorical structures are ubiquitous in various areas
of mathematics, including algebraic topology, representation theory and
rewriting theory. One of the most general setups is the one of higher
categories which consist of objects and morphisms in each dimension with
composition laws along lower-dimensional morphisms that satisfy some version
of the usual axioms. Depending on whether the axioms are satisfied on the
nose or up to coherent higher invertible morphisms, one talks about
\emph{strict} or \emph{weak} $\omega$-categories, the latter
ones being also known as $(\infty,\infty)$-categories.
The definition for strict $\omega$-categories
is well established (see for instance \cite{StreetOrientedSimplexes}). On
the contrary, many approaches for weak $\omega$-categories compete (see for
instance~\cite{BataninWCat, PenonApproche, StreetWeakOmega,futia,Lumsdaine}) and
are not clearly compared yet.

Although homotopy theorists are ultimately interested in weak
$\omega$-categ\-ories, strict $\omega$-categories have their importance in
homotopy theory for several reasons. First, as shown by the second author
\cite{GagnaStrict}, they provide a model for homotopy types, which has been
studied extensively by Maltsiniotis and the first author, for instance in
\cite{AraMaltsiniotisVers, AraMaltsiniotisThmASimp,
AraMaltsiniotisThmACat}. Second, they give an easier setting to
develop or test ideas for the general theory of higher categories, and
provide a source of intuition. Third, and most importantly for homotopy
theorists, they contain diagrammatic shapes encoding the various ways of
composing a bunch of cells, leading to several definitions of weak
$\omega$-categories. The idea of basing a definition for weak
$\omega$-categories on strict shapes dates back to Grothendieck's definition
of weak $\omega$-groupoids \cite{GrothPS}, and was used for various definitions of weak
$\omega$-categories \cite{BataninWCat, PenonApproche, MaltsiniotisGroupoids}. Similarly, Joyal's
category $\Theta$ \cite{JoyalDisks}, which plays a key role in higher category theory,
can be seen as a category of strict $\omega$\nobreakdash-categories.
Although the shapes of these $\omega$-categories are pretty
simple, other similar definitions are based on more complicated shapes, such
as in the work of Street and Verity \cite{StreetOrientedSimplexes, VerityComplicialAMS}, Al-Agl, Brown and Steiner \cite{AlAglBrownSteiner}, Barwick and Schommer-Pries \cite{BarwickSchommerPries}, Henry \cite{HenryRegular} or Hadzihasanovic
\cite{HadzihasanovicDiagrammatic}.

Like for ordinary algebra, an important role amongst all strict $\omega$-categories is played by those 
that are freely generated by gluing a number of cells of different dimensions. The datum describing such a free strict $\omega$-category was introduced independently by Burroni \cite{BurroniWordProblems} under the name of \emph{polygraph} and by Street in the $2$-categorical context \cite{StreetOrientedSimplexes} under the name of \emph{computad}, with the goal of dealing with manageable presentations of higher categories.

Polygraphs are the $\omega$-categorical analog of CW-complexes in topology, being built out of cells glued along their boundaries. This is formalized by the fact that $\omega$-categories generated by polygraphs are precisely the cellular objects in the canonical model structure of strict $\omega$-categories by Lafont--M\'etayer--Worytkiewicz \cite{LMW}. Better, Métayer proved that they are exactly the cofibrant objects of this model structure~\cite{MetCof}. 
They also contain all the cofibrant objects of the conjectured Thomason model structure for homotopy types by Maltsiniotis and the first author \cite{AraMaltsiniotisVers}. Furthermore, polygraphs provide the resolutions to compute the so-called polygraphic homology of $\omega$-categories \cite{LafontMetayer}, a construction that generalizes classical homology of groups, monoids or categories \cite{LafontMetayer, GuettaHomology}.

Polygraphs also play an important role in rewriting theory.
While we know that, in general, the word problem in an arbitrary group is not soluble, in some cases it is possible to get (more or less efficient) algorithms for manipulating words in a given generating system by orienting the relations.
In order to figure out the relation between different relations and how coherent the whole rewriting procedure is, a very successful formalization is given by polygraphs.
In particular, polygraphic homology is closely related to homology computations via Squier theory for rewriting systems. For a modern treatment of Squier theory, we refer the reader e.g.~to~\cite{GuiraudMalbosHigherdimensional,GuiraudMalbosPolygraphs}.

 It is hard to encode in a concise way the information of an $n$-category, and even more so of an $\omega$-category.
A manifestation of this complication is that even understanding maps between most fundamental, elementary objects is oftentimes combinatorially highly involved.
 Steiner \cite{SteinerEmbedding} identified a class of particularly nice  $\omega$\nobreakdash-categories --- which we will refer to as \emph{strong Steiner $\omega$\nobreakdash-categories} --- that admit an algebraic model, in the form of a chain complex with extra structure.
 
 For this class of strict $\omega$-categories, significant progress was made using Steiner's machinery. In the context of strict $\omega$-categories, the second author used strong Steiner $\omega$-categories to perform key computations in establishing that strict $\omega$-categories model homotopy types \cite{GagnaStrict}, and Maltsiniotis and the first author studied the Gray tensor product and developed 
 a theory of join and slices
for strict $\omega$-categories \cite{AraMaltsiniotisJoin}, and proved versions of Quillen's Theorems A and B~\cite{AraMaltsiniotisThmASimp,AraThmB,AraMaltsiniotisThmACat}.
 In the context of (weak) complicial sets, Steiner~\cite{SteinerOrientals,SteinerUniversal}
 gave algebraic descriptions of orientals and maps between them, Medina-Mardones \cite{MedinaMardonesGlobular} recognized the connection between orientals and Steenrod's $\cup_i$-products, Loubaton realized strict $\omega$-categories as weak complicial sets \cite{LoubatonNerfs}, Maehara showed that orientals provide fibrant replacements of simplices \cite{MaeharaOrientals}, and the last three authors showed that the nerve construction is compatible with two flavors of join constructions~\cite{GOR1}.

For a strict $\omega$-category $\cC$, being a strong Steiner $\omega$-category relies on a freeness condition (recalled as \cref{defnbasis}) and a loop-freeness condition (recalled as \cref{defn:Steiner_loopfree_wcats}) that are both expressed in terms of the chain complex~$\lambda\cC$ obtained by ``linearizing'' $\cC$. Steiner proved in \cite[Theorem~6.1]{SteinerEmbedding} that strong Steiner $\omega$-categories are freely generated by polygraphs. Nevertheless, while writing the papers \cite{GOR1} and \cite{ALM}, we realized that it was not clear how Steiner's notion of freeness relates to freeness in the sense of polygraphs, and how Steiner's notion of loop-freeness could be expressed without explicit reference to the chain complex $\lambda\cC$.

Our goal in this paper is to characterize this class of $\omega$-categories in purely categorical terms.
 Although the final description is not utterly surprising, the relationship between the linear and $\omega$-categorical viewpoint is subtle and tricky, and we hope to ``demystify''
 this aspect with this work.
 Our main result is that strong Steiner $\omega$-categories are exactly those $\omega$-categories freely generated by a polygraph that in addition satisfy a loop-freeness condition that we introduce (as \cref{stronglyloopfreepolygraph}) described purely in terms of the categorical structure:

\begin{unnumberedtheorem}
For an $\omega$-category $\cC$, the following are equivalent.
\begin{enumerate}[leftmargin=*]
    \item The $\omega$-category $\cC$ is freely generated by a strongly loop-free polygraphic generating set.
    \item  The $\omega$-category $\cC$ is a strong Steiner $\omega$-category, i.e., it admits a strongly loop-free atomic basis in the sense of \cite{SteinerEmbedding}.
\end{enumerate}
\end{unnumberedtheorem}

This theorem, which is part of \cref{equivalencestrongSteiner}, gives a very satisfactory 
characterization of strong Steiner $\omega$\nobreakdash-categories giving an easier criterion to check whether or not a given $\omega$\nobreakdash-categ\-ory is a strong Steiner $\omega$-category.

Our paper is organized as follows: in \cref{sec:Recollections} we recall the necessary background of Steiner theory, and in \cref{sec:Categorical} we define categorical versions of the involved notions, leading up to the desired characterization.

 \subsection*{Acknowledgements}
 We are thankful to L\'eonard Guetta for helpful conversations.
The second-named author gratefully acknowledges the support of Praemium Academiae of M.~Markl and RVO:67985840. The authors would like to thank the MFO, Oberwolfach, and the CIRM, Marseille, for the excellent environments in which part of the material of this paper was developed.
We are thankful to the referee for useful comments.

 \section{Recollection of Steiner's theory}
 \label{sec:Recollections}
 
 In this section we recall the main ingredients of Steiner's theory from \cite{SteinerEmbedding}, also later recalled in \cite[\textsection2]{AraMaltsiniotisJoin}. There is a class of nice $\omega$-categories --- referred to as \emph{strong Steiner $\omega$-categories} --- and a nice class of augmented directed chain complexes --- referred to as \emph{strong Steiner complexes} --- that fit into an equivalence of categories via the adjoint pair $(\lambda,\nu)$.
 
\subsection{\pdfomega-categories}

While we refer the reader to e.g.~\cite{StreetOrientedSimplexes} for a traditional approach to the definition of an $\omega$-category, we briefly recall the main features here.

The data of an \emph{$\omega$-category} $\cC$ consists of a collection of sets $\cC_q$, for ${q \geq 0}$,
where $\cC_0$ is called the set of \emph{objects} of $\cC$ and $\cC_q$ for $q>0$
is the set of \emph{$q$-cells} or cells of \emph{dimension} $q$ of $\cC$, together with:
\begin{itemize}[leftmargin=*]
    \item \emph{source} and \emph{target} operators $s_q, t_q \colon \cC_p \to \cC_q$
    for all $p \geq q \geq 0$;
    \item \emph{identity} operators $\id_q \colon \cC_p \to \cC_{q}$ for all $q\geq p\ge0$;
    \item \emph{composition} operators $\comp_p \colon \cC_q \times_{\cC_p} \cC_q \to \cC_q$
    defined for all $q > p \geq 0$ and all pairs of $q$-cells $(g, f)$ for which $s_p(g) = t_p(f)$.
\end{itemize}
We say that $\cC$ is an $\omega$-category if for all $r > q > p \geq 0$ the triple $(\cC_p, \cC_q, \cC_r)$
together with all the relevant source, target, identity and composition operators
is a $2$-category. In particular, 
\begin{equation}
\label{globularity}
s_ps_q(f)=s_p(f)    
\qquad\text{and}\qquad
t_pt_q(f)=t_p(f)    
\end{equation}
for any $r$-cell $f$ of $\cC$ and $r>q>p$.

An \emph{$\omega$-functor} $F \colon \cC \to \cD$ between $\omega$-categories $\cC$ and $\cD$ is a collection of maps $F_q \colon \cC_q \to \cD_q$ for
$q \geq 0$ that preserves source, target, identity, and composition operators.
We denote by $\omegacat$ the category of (small) $\omega$\nobreakdash-categories and $\omega$-functors.

A cell in an $\omega$-category $\cC$ is said to be \emph{trivial} if it is the identity of a cell of lower dimension. For $n\ge0$, an \emph{$n$-category} is an $\omega$-category in which all $q$-cells are trivial for $q>n$, and an \emph{$n$-functor} is an $\omega$-functor between $n$-categories. We denote by $n\cat$ the (full) subcategory of $\omega\cat$ given by $n$-categories and $n$-functors.

\begin{rmk}
\label{truncation}
For $n>0$, the canonical inclusion $n\cat\hookrightarrow \omega\cat$ admits a right adjoint $\varpi_{n}\colon \omega\cat\rightarrow n\cat$, which produces an $n$-category $\varpi_{n}\cC$ by forgetting all non-trivial $q$-cells of $\cC$ for $q> n$
and does not change the underlying $n$-category of $\cC$. This functor is called the \emph{$n$-truncation} functor and treated e.g.~in~\cite[\textsection 1.2]{AraMaltsiniotisJoin}.
\end{rmk}

    The $n$-cells of an $\omega$-category are corepresented by
    the $n$-disk $\omega$-category:
    
\begin{ex}
\label{ex:disks}
For $n\ge0$, the \emph{$n$-disk} $\cD[n]$  is the $n$-category having exactly one non-trivial
    $n$-cell $\varsigma$ and two distinct non-trivial $q$-cells $s_q(\varsigma)$ and $t_q(\varsigma)$
    for every $0 \leq q < m$.
    For the first few values of $n$, the corresponding disk can be depicted as follows\footnote{In this example and later, when drawing an $\omega$-category, we follow the convention that we only draw the generating non-trivial cells, as arrows of the correct dimension pointing from their source to their target.}
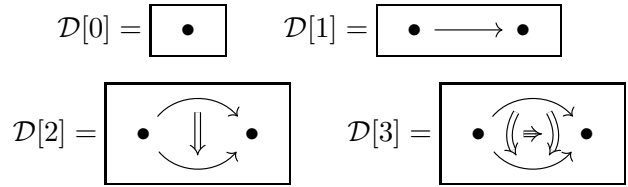
    \[
			\cD[0] =
			\fbox{
			\begin{tikzcd}[ampersand replacement=\&]
				\bullet
			\end{tikzcd}
		    }
		\quad\quad
			\cD[1] =
			\fbox{
			\begin{tikzcd}[ampersand replacement=\&]
				\bullet \arrow[r] \& \bullet
			\end{tikzcd}
			}
			\ \ 
			 		\]
		\[\cD[2] =
			\fbox{
\begin{tikzcd}[ampersand replacement=\&]
				\bullet \arrow[r, bend left=50, ""{name=U, below}] \arrow[r, bend right=50, ""{name=D}]
				\& \bullet
				\arrow[Rightarrow,from=U,to=D]
			\end{tikzcd}
			}
		\quad\quad
			\cD[3] = 
			\fbox{
			\begin{tikzcd}[ampersand replacement=\&]
				\bullet \arrow[r, bend left=50, ""'{name=U}] \arrow[r, bend right=50, ""{name=D}]
				\& \bullet
				\arrow[Rightarrow,from=U,to=D, shift right=1ex, bend right=30, ""{name=L}]
				\arrow[Rightarrow, from=U, to=D, shift left=1ex, bend left=30, ""{name=R, left}]
				\arrow[triple, from=L, to=R]{}
			\end{tikzcd}
			}
		\]
\end{ex}\ 

Similarly, parallel $n$-morphisms are corepresented by the $n$-sphere:

\begin{ex}
\label{ex:sphere}
For $n\ge0$, the \emph{$(n-1)$-sphere} $\cS[n-1]$ is the $(n-1)$-truncation $\varpi_{n-1}\cD[n]$ of the $n$-disk, that is, the $(n-1)$-category having exactly two distinct non-trivial $q$-cells $s_q(\varsigma)$ and $t_q(\varsigma)$
    for every $0 \leq q < n$.
    For the first few values of $n$, the corresponding sphere can be depicted as follows
    \[
			\cS[-1] = \fbox{
		\vphantom
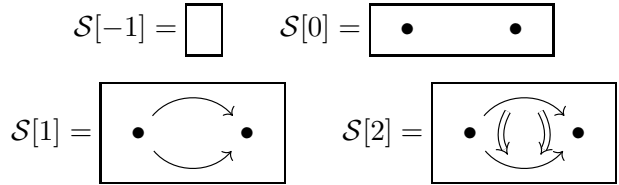
{\begin{tikzcd}[ampersand replacement=\&]
				\bullet  \& \bullet
			\end{tikzcd}}
			}
			\quad\quad 
			\cS[0] =
			\fbox{
			\begin{tikzcd}[ampersand replacement=\&]
				\bullet  \& \bullet
			\end{tikzcd}
			}
	\]
	\[
			\cS[1] = 
			\fbox{
			\begin{tikzcd}[ampersand replacement=\&]
				\bullet \arrow[r, bend left=50, ""{name=U, below}] \arrow[r, bend right=50, ""{name=D}]
				\& \bullet
			\end{tikzcd}
			}
\quad\quad
			\cS[2] = \fbox{
			\begin{tikzcd}[ampersand replacement=\&]
				\bullet \arrow[r, bend left=50, ""'{name=U}] \arrow[r, bend right=50, ""{name=D}]
			\& \bullet
				\arrow[Rightarrow,from=U,to=D, shift right=1ex, bend right=30, ""{name=L}]
				\arrow[Rightarrow, from=U, to=D, shift left=1ex, bend left=30, ""{name=R, left}]
			\end{tikzcd}\
			}
		\]
\end{ex}\ %

An important family of examples of $n$-categories is the collection of objects of Joyal's cell category $\Theta_n$ \cite{JoyalDisks}.

\begin{ex}
\label{ex:theta}
For $n\ge0$, the objects of Joyal's category $\Theta_n$ are $n$-categories describing some ``globular pasting schemess''. Let us be more specific in low dimension.
The only object of the category $\Theta_0$ is the terminal $\omega$-category
\[
    [0] =
	\fbox{
	\begin{tikzcd}[ampersand replacement=\&]
		\bullet
	\end{tikzcd}
    }
\]
and $\Theta_1$ is the simplex category $\Delta$ seen as a full subcategory of $\cat$, whose generic object is the ordinal category $[m]$.
\[
	[m] =
	\fbox{
	\begin{tikzcd}[ampersand replacement=\&]
		\bullet \arrow[r] \& \bullet \arrow[r] \& \cdots \arrow[r] \& \bullet
	\end{tikzcd}
	}
\]
The objects of the category $\Theta_2$ are of the form
$[m|k_1,\dots,k_m]$, for $m\ge0$ and $k_i\ge0$, where for instance the
$2$-category $[3|2,0,1]$ is (generated by) the following data.
   \[
[3|2,0,1]=
\fbox{\begin{tikzcd}[row sep=3.2cm, column sep=2.2cm,ampersand replacement=\&]
  \bullet \arrow[r, bend left=50, "", ""{name=U,inner sep=2pt,below}]
  \arrow[r, ""{near end, xshift=0.2cm}, ""{name=D,inner sep=2pt},""{name=M,inner sep=2pt, below}]
  \arrow[r, bend right=50, ""{below}, ""{name=DD,inner sep=2pt}]
  \& \bullet\arrow[r, ""]
    \& [-0.9cm]
    \bullet\arrow[r, bend left=50, "" ""{name=U1,inner sep=2pt,below}]
  \arrow[r, ""{near end, xshift=0.2cm}, ""{name=D1,inner sep=2pt},""{name=M1,inner sep=2pt, below}]
  \& \bullet
    \arrow[Rightarrow, from=U, to=D, ]
    \arrow[Rightarrow, from=M, to=DD, ""]  
  \arrow[Rightarrow, from=U1, to=D1, ""]
\end{tikzcd}
}\]
\end{ex}\ %

Finally, another important family of examples of $n$-categories is given by
Street's orientals \cite{StreetOrientedSimplexes}.

\begin{ex}
\label{ex:oriental}
For $n\ge0$, the $n$-th oriental $\cO[n]$ is an $n$-category somehow shaped as an $n$-simplex.
For the first few values of $n$, the corresponding oriental can be depicted as follows.
      \[
			\cO[0] =
			\fbox{
			\begin{tikzcd}[ampersand replacement=\&]
				\bullet
			\end{tikzcd}
			}
\quad\quad
			\cO[1] =
		\fbox{	\begin{tikzcd}[ampersand replacement=\&]
				\bullet \arrow[r] \& \bullet
			\end{tikzcd}
			}
		\quad\quad
		   	\cO[2] = 
				\fbox{
			 \begin{tikzcd}[row sep=small, column sep=tiny,ampersand replacement=\&]
                                            \& \bullet \ar[rd]   \&         \\
        \bullet \ar[ru] \ar[rr, ""{name=s}] \&                   \& \bullet
                \ar[Rightarrow, from=s, to=1-2, shorten >= 2pt]
        \end{tikzcd}
}
		\]
\begin{center}
    \centering
    \begin{tikzpicture}
    \draw (-3.0,-1) rectangle (3.0,1);
			\square{
				/square/label/.cd,
	     			0=$\bullet$, 1=$\bullet$, 2=$\bullet$, 3=$\bullet$,
	     			01={}, 12={}, 23={}, 03={},
	     			012={}, 023={}, 123={}, 013={},
	     			0123={}
     			}
	\end{tikzpicture}
\end{center}
For a precise account, we refer the reader to Street's original construction \cite{StreetOrientedSimplexes} of~$\cO[n]$ or to Steiner's alternative construction~\cite{SteinerOrientals}.
\end{ex}

\label{sec:Steiner}

\subsection{Augmented directed chain complexes}
By a \emph{chain complex} $C$ we will always mean an
$\mathbb{N}$-graded chain complex of abelian
groups with homological indexing, that is,
a family $(C_q)_{q\geq 0}$ of abelian groups, together with maps $\partial_{q} \colon C_{q+1} \to C_{q}$ satisfying
\[\partial_q \partial_{q+1}=0.\]

Given chain complexes $C$ and $\overline C$, a \emph{chain map} or \emph{morphism of chain complexes} $\phi\colon C\to \overline C$ consists of a family of homomorphisms $(\phi_q\colon C_q\to \overline C_q)_{q\geq 0}$ that commutes with the differentials in the sense that
\[\overline\partial_{q} \phi_{q+1}=\phi_q\partial_{q}\]
for every $q \geq 0$.

An \emph{augmented chain complex} is a pair $(C,\varepsilon)$
of a chain complex $C$ and an augmentation, namely a map $\varepsilon\colon C_0\to\mathbb Z$ such that \[\varepsilon\partial_0=0.\]

An \emph{augmented chain map} $\phi \colon (C,\varepsilon) \to (\overline C, \overline\varepsilon)$ between augmented chain complexes $(C,\varepsilon)$ and $(\overline C,\overline\varepsilon)$ consists of a chain map $\phi\colon C\to\overline C$ that is moreover compatible with the augmentations, namely such that
\[\overline\varepsilon \phi_0 = \varepsilon.\]

Steiner \cite[\textsection2]{SteinerEmbedding} introduced an enhancement of the structure of an augmented chain complex:

\begin{defn}[{\cite[Def.~2.2]{SteinerEmbedding}}]
 An \emph{augmented directed complex}
 is a triple
 $(C, C^+, \varepsilon)$ where $(C, \varepsilon)$ is
 an augmented chain complex and $C^+ =  (C^+_q)_{q\ge0}$ is a collection of commutative monoids,
 with $C^+_q$ is a submonoid of $C_q$ called the
 \emph{positivity submonoid} of $C_q$.
 
 Given augmented directed chain complexes
$(C, C^+, \varepsilon)$ and $(\overline C, \overline C^+, \overline\varepsilon)$,
 a \emph{morphism of augmented directed chain complexes} or \emph{augmented directed chain map} $\phi \colon (C, C^+, \varepsilon) \to (\overline C, \overline C^+, \overline\varepsilon)$ is an augmented chain map $\phi \colon (C, \varepsilon) \to (\overline C, \overline{\varepsilon})$ that moreover preserves the positivity submonoids\footnote{The differentials of an augmented directed chain complex
    \emph{need not} respect the positivity submonoid.}, namely satisfies
    \[\phi_q(C^+_q) \subseteq \overline C^+_q\]
 for all $q\ge0$.
\end{defn}
We denote by $\adch$ the category of augmented directed chain complexes and augmented directed chain maps.

\subsection{Steiner's functors}
\label{sec:nu}

Steiner \cite[\textsection2]{SteinerEmbedding} constructed a pair of adjoint functors
\[\lambda\colon \omega\cat\rightleftarrows\adch\cocolon\nu.\]
This means that for any $\omega$-category $\cC$ and any augmented directed chain complex $C'$ there is a natural bijection
\[\adch(\lambda \cC,C')\cong\omega\cat(\cC,\nu C').\]
To give a bit of context, we recall the basic data of these constructions, borrowing some notations from \cite[\textsection2.4]{AraMaltsiniotisJoin}, and refer the reader to the original sources for more details.

For $C$ an augmented directed chain complex, 
the set of $q$-cells $(\nu C)_q$ of the
$\omega$\nobreakdash-category $\nu C$ for $q \geq 0$ is the set of
tables
\[
	x=\tabld{x}{q}
\]
such that, for $\alpha=+,-$ and $0\le p\le q$, the following hold:
\begin{enumerate}
	\item $x^\alpha_p$ belongs to $C^+_p$;
	\item $\partial_{p-1}(x^\alpha_p) = x^+_{p-1} - x^-_{p-1}$ for $0<p\leq q$;
	\item $\varepsilon(x^\alpha_0) = 1$;
	\item $x_q^- = x_q^+$.
\end{enumerate}
We refer to~\cite[Def.~2.8]{SteinerEmbedding} or to~\cite[\textsection 2.4]{AraMaltsiniotisJoin}
for a full description of the $\omega$-categorical structure of $\nu C$.

For $\cC$ an $\omega$-category, we recall the definition of the augmented directed chain complex $\lambda\cC$ from \cite[Def.~2.4]{SteinerEmbedding}. For $q\ge0$, the abelian group of $q$-chains of $\lambda\cC$ is the quotient of $\mathbb Z[\cC_q]$ given by
\begin{equation}
\label{lambda}
(\lambda\cC)_q \coloneqq \frac{\mathbb Z[\cC_q]}{\left<[x\ast_p y]_q-[x]_q-[y]_q\mid x,y\in\cC_q;p<q\right>},
\end{equation}
where, for $z$ a $q$-cell of $\cC$, we denoted by $[z]_q$ the corresponding element of~$\mathbb{Z}[\cC_q]$.
The positivity submonoid $(\lambda\cC)_q^+$ is the submonoid of $(\lambda\cC)_q$ generated by the collection of elements $[f]_q$ for $f$ a $q$-cell of $\cC$.
The boundary maps $\partial_{q-1}\colon(\lambda\cC)_q\to(\lambda\cC)_{q-1}$ are determined by the condition on generators
\[\partial_{q-1}([f]_{q}) \coloneqq [t_{q-1}(f)]_{q-1}-[s_{q-1}(f)]_{q-1},\]
where $f$ is in $\cC_q$, and the augmentation map $\varepsilon\colon(\lambda\cC)_0\to\mathbb Z$ by the condition on generators
\[\varepsilon([x]_0) \coloneqq 1,\]
where $x$ is in $\cC_0$.

Let's illustrate a couple of general phenomena about the construction $\lambda\cC$ by exploring an explicit example.

\begin{ex}
\label{examplelambda}
Consider the $\omega$-category $\cO[3]$.
\begin{center}
    \centering
    \begin{tikzpicture}
    \draw (-3.0,-1) rectangle (3.0,1);
			\square{
				/square/label/.cd,
	     			0 = {$x$}, 1={$y$}, 2={$z$}, 3={$w$},
                    01={$f$}, 12={$g$}, 23={$h$},
                    02={$i$}, 13={$j$}, 03={$k$},
                    023={$\alpha$}, 012={$\beta$},
                    013={$\gamma$}, 123={$\delta$},
                    0123={$\Gamma$}
     			}
	\end{tikzpicture}
\end{center}
\begin{enumerate}[leftmargin=*]
    \item Although $[f]_1\neq0\in(\lambda\cO[3])_1$, we have
    $[\id_2(f)]_2=0\in(\lambda\cO[3])_2$.
    \item We have $[\beta]_2=[\id_2(h) \ast_0 \beta]_2$.
    \item We have $[(\id_2(h) \ast_0 \beta) \ast_1 \alpha]_2=\cancel{[\id_2(h)]_2}+[\beta]_2+[\alpha]_2=[\beta]_2+[\alpha]_2$.
\end{enumerate}
\end{ex}

The facts discussed in the example are instances of general facts.

\begin{rmk}
\label{propertieslambda}
Let $\cC$ be an $\omega$-category.
\begin{enumerate}[leftmargin=*]
    \item If a $q$-cell $x$ in $\cC$ is trivial then the class $[x]_q$ vanishes in $(\lambda\cC)_q$, as shown in \cite[Proposition~2.5]{SteinerEmbedding}.
    \item If two $q$-cells of $\cC$ differ by a $(q-1)$-dimensional whiskering, they represent the same class in $(\lambda\cC)_q$.
    \item Roughly speaking, in $\lambda\cC$ composition becomes addition.
\end{enumerate}
\end{rmk}

The functors $\lambda$ and $\nu$ do not define equivalences of categories in the treated generality but we'll see in \cref{SteinersThm} that they do induce equivalences of categories when restricted and corestricted to suitable full subcategories.

\subsection{Strong Steiner complexes}
\label{sec:steiner_complexes}
Amongst all augmented directed chain complexes, Steiner \cite{SteinerEmbedding} identified a class of particularly nice ones, which we refer to as \emph{strong Steiner complexes}, adopting the terminology from \cite[\textsection2]{AraMaltsiniotisJoin}. This notion builds on a few preliminary definitions and notation, which we briefly recall.

\begin{defn}[{\cite[Definition~3.1]{SteinerEmbedding}}]
Let $(C, C^+, \varepsilon)$ be an augmented directed chain complex, and $B = (B_q)_{q\ge0}$ a family of subsets $B_q\subseteq C_q$ for~$q\ge0$. The $\mathbb{N}$-graded set $B$ is a \emph{basis} for $C$ if for every $q \geq 0$ the set $B_q \subseteq C_q$ is
\begin{itemize}[leftmargin=*]
    \item  a basis of the abelian group $C_q$, namely $C_q\cong\mathbb Z[B_q]$, and
    \item a basis of the commutative monoid $C^+_q$, namely $C_q^+\cong\mathbb N[B_q]$.
\end{itemize} 
\end{defn}

If $C$ has a basis $B$, then there are inclusions $B_q\subseteq C^+_q\subseteq C_q$.

\begin{lem}[{\cite[\textsection3]{SteinerEmbedding}}]
\label{basiscomplexunique}
    If an augmented directed chain complex $C$ admits a basis $B$,
    then the basis is uniquely determined.
\end{lem}
    
    In particular, if $C$ admits a basis, it makes sense to talk about \emph{the} basis of $C$.
    
    \medbreak

The following was introduced in \cite[\textsection2]{SteinerEmbedding}, and also treated in \cite[\textsection2.7]{AraMaltsiniotisJoin}.

\begin{defn}
\label{SupportComplex}

    Let $(C, C^+, \varepsilon)$ be an augmented directed chain complex
    with basis $B$. For $q \geq 0$, an element $c$ in $C_q$
    has a canonical decomposition
    \begin{equation}
\label{canonicaldecomposition}
        c = \sum_{b \in B_q} \lambda_b \cdot b,
        \quad \lambda_b \in \mathbb{Z}.
    \end{equation}
    The \emph{positive support} and \emph{negative support} of $c$ are the finite subsets of $B_q$ given respectively by
    \[\supp^+c=\{ b \mid \lambda_b>0\}\quad\text{and}\quad\supp^-c=\{ b \mid \lambda_b<0\}.\]
    The \emph{support} $\supp(c)$ is the union of the positive and negative support, that is, the set of generators appearing in the linear expansion of $c$ with respect to the basis~$B_q$. The \emph{positive part} and
    \emph{negative part} of $c$ are respectively
    \[
        c^+ \coloneqq \sum_{\substack{b \in B_q \\ b \in\supp^+c}} \lambda_b \cdot b
       \in C_q^+ \quad\text{and}\quad
        c^- \coloneqq -\sum_{\substack{b \in B_q \\ b \in\supp^-c}} \lambda_b \cdot b\in C_q^+.
    \]
    We denote
    \[\partial_{q-1}^+c \coloneqq (\partial_{q-1} c)^+\quad\text{and}\quad\partial_{q-1}^-c \coloneqq (\partial_{q-1} c)^-.\]
\end{defn}

The following symbol keeps track of the iterated positive and negative parts of the boundaries of a chain.

\begin{notn}[{\cite[\textsection2.8]{AraMaltsiniotisJoin}}]
\label{def:atom}
	Let $(C, C^+, \varepsilon)$ be an augmented directed chain complex.
	For $\alpha=+,-$ and $c\in C_q$, the symbol $\atom{c}^{\alpha}_{p}$ is defined inductively on $p=q,\dots,0$ by
	\[\atom{c}^{\alpha}_p\coloneqq\left\{\begin{array}{cl}
	    c & \text{if }p=q\  \\
	    \partial_{p}^\alpha(\atom{c}^\alpha_{p+1}) & \text{if }0\le p<q.
	\end{array}\right.\]
\end{notn}

	\begin{defn}[{\cite[Definition~3.4]{SteinerEmbedding}}]
	\label{defn:unital}
	Let $(C, C^+, \varepsilon)$ be an augmented directed chain complex with basis $B = (B_q)_{q\ge 0}$.
	The basis $B$ of $C$ is said to be \emph{unital} if
	\[\varepsilon(\atom{c}^-_0) = 1 = \varepsilon(\atom{c}^+_0)\]
	for all $c$ in $B_q$ and $q\ge0$.
\end{defn}

\begin{notn}[{\cite[Definition~3.6]{SteinerEmbedding}}]
\label{loopfree_preorder_complexes}
 Let $C$ be an augmented directed chain complex with basis $B$. We denote by $\leq_{\mathbb N}$ the preorder relation on $\coprod_{q\ge0} B_q$ generated by the condition
        \[\text{$a\leq_{\mathbb N}b$ if $a\in B_p$, $b\in B_q$, and
        $\begin{cases}\text{either $a\in \supp(\partial^-_{q-1}(b))$ with $q>0$,}\\
        \text{or $b\in \supp (\partial^+_{p-1}(a))$ with $p>0$.}\end{cases}$} \]
\end{notn}

\begin{defn}[{\cite[\textsection3.6]{SteinerEmbedding}}]
\label{defbasisalgebraic}
    Let $C$ be an augmented directed chain complex with a unital basis. If the preorder relation $\leq_{\mathbb{N}}$ is a partial order, we say that $C$ admits a \emph{strongly loop-free unital basis}.
\end{defn}

The following notion plays a key role in Steiner's theory \cite[\textsection\textsection4-5]{SteinerEmbedding} and hence was given a name in~\cite[Ch.~2]{AraMaltsiniotisJoin}. 

\begin{defn}[{\cite[\textsection2.15]{AraMaltsiniotisJoin}}]
    An augmented directed chain complex with a strongly loop-free unital basis is a \emph{strong Steiner complex}.
\end{defn}

\subsection{Strong Steiner \pdfomega-categories}

Steiner \cite{SteinerEmbedding} also identified which $\omega$\nobreakdash-categories correspond, in a precise sense, to strong Steiner complexes via the adjunction $(\lambda,\nu)$. Following \cite[\textsection2]{AraMaltsiniotisJoin}, we'll refer to these $\omega$-categories as \emph{strong Steiner $\omega$-categories}.

\begin{defn}[{\cite[Definition~4.1]{SteinerEmbedding}}]
    Let $\cC$ be an $\omega$-category, and $E$ a set of cells in $\cC$. We say that the subset $E$ of the cells of $\cC$ \emph{composition-generates} $\cC$ if the smallest subset containing $E$ and closed under composites and identities is the set of all cells of $\cC$.
\end{defn}

In particular, if $E$ composition-generates $\cC$, then $E$ must contain all the objects of $\cC$
    and every cell of $\cC$ can be written as a finite composition of cells of $E$ or
    iterated identities of cells of $E$.

\begin{defn}[{\cite[Definition~4.4]{SteinerEmbedding}}]
\label{defnbasis}
    Let $\cC$ be an $\omega$-category, and $E$ a subset of the cells of $\cC$
    which composition-generates $\cC$. We say that the subset $E$ of the cells of $\cC$
    is a \emph{basis} for $\cC$ if for every $q\ge 0$:
    \begin{itemize}[leftmargin=*]
    \item the assignment $e\mapsto[e]_q$ defines a bijection between $E\cap\cC_q$ and the set $[E\cap\cC_q] \coloneqq \{[e]_q\in(\lambda\cC)_q\mid e\in E\cap\cC\}$; 
        \item the set $[E\cap\cC_q]$ is a basis of the abelian group $(\lambda \cC)_q$, namely we have $(\lambda \cC)_q\cong\mathbb Z[E\cap\cC_q]$; and
        \item the set $[E\cap\cC_q]$ is a basis of the monoid $(\lambda \cC)^+_q$, namely $(\lambda \cC)_q^+\cong\mathbb N[E\cap\cC_q]$.
    \end{itemize}
\end{defn}

\begin{rmk}
\label{positivecoefficients}
If $E$ is a basis for $\cC$, then for any $q$-cell $f$ of $\cC$ the coefficients $\lambda_b\in\mathbb Z$ occurring in the decomposition \eqref{canonicaldecomposition} of $[f]_q$ into elements of the basis $[E]$,
\[[f]_q=\sum_{b \in E\cap\cC_q} \lambda_b \cdot [b]_q,\]
are all non-negative, i.e., we have $\lambda_b\ge0$.
\end{rmk}

The following proposition allows us to speak of \emph{the} basis of
an $\omega$-category.

\begin{prop}
\label{UniqueOmegaBasis}
  If an $\omega$-category admits a basis $E$,
    then $E$ is uniquely determined.
\end{prop}

\begin{proof}
Let $\cC$ be an $\omega$-category and let $E$, $E'$ be two bases for $\cC$. In particular we know that $[E]$ and $[E']$ are bases of $\lambda\cC$, and are therefore equal by \cref{basiscomplexunique}.
In particular, from the definition of basis we obtain a bijection
$E\cap\cC_q\cong[E\cap\cC_q]=[E'\cap\cC_q]\cong E'\cap\cC_q$.

In order to prove that $E'\subseteq E$,
we assume for contradiction that $E'\nsubseteq E$. Let $a'$ be a cell in $E'\setminus E$ of minimal dimension $q\ge0$, and let $a$ be the unique $q$-cell in $E$ corresponding to $a'$ via the bijection $E\cap\cC_q\cong E'\cap\cC_q$. Namely, $a$ is the unique $q$-cell in $E$ such that $[a]_q=[a']_q$. Since both $E$ and~$E'$ composition-generate $\cC$, all objects of $\cC$ must belong to $E$ and $E'$, and so $q>0$. Given that $E$ composition-generates $\cC$, there is an expression for $a'$ as a composite of trivial $q$-cells in $E$ and exactly $k$ occurrences in total of non-trivial (not necessarily distinct) $q$-cells in $E$ for some $k\ge0$, say  $a_1,\dots,a_k$. We then obtain in $(\lambda\cC)_q$ two expressions of $[a']_q$ in terms of elements of $[E\cap\cC_q]$:
\[[a]_q=[a']_q=[a_1]_q+\dots+[a_k]_q.\]
Since $[E\cap\cC_q]$ is a basis of $(\lambda\cC)_q$, the expression \eqref{canonicaldecomposition} for $[a']_q$ in terms of elements of $[E\cap\cC_q]$ is unique, and the only possibility is that $k=1$ and~$a_1=a$. This means that $a'$ is a composite of (exactly one occurrence of) $a$ with other trivial $q$-cells in $E$. With a symmetric argument, we also deduce that $a$ is a composite of (exactly one occurrence of) $a'$ with other trivial $q$-cells in $E'$.

We can thus substitute the expression for $a'$ into the expression for $a$ and obtain an expression for $a$ as a composite of (exactly one occurrence of) $a$ and other trivial $q$-cells of $\cC$. 
Let $0<p<q$ be maximal such that $c$ is a non-trivial $p$-cell for which $\id_{q}(c)$
occurs in the expression of $a$ mentioned above. By taking $s_p$, we obtain an expression of $s_p(a)$ as a composite involving at least one occurrence of $s_p(a)$, at least one occurrence of $c$ and possibly other $p$-cells.
This leads to an expression of~$[s_p(a)]_p$
as a positive integer linear combination involving at least one occurrence of $[s_p(a)]_p$, at least one occurrence of $[c]_p$ and possibly the $p$-classes of other $p$-cells. By canceling~$[s_p(a)]$ on both sides, and using \cref{positivecoefficients}, we obtain
 a positive non-trivial linear combination of elements in $[E\cap\cC_p]$ which vanishes in $ (\lambda \cC)_p$, contradicting the fact that $[E\cap\cC_p]$ is a basis for the abelian group $(\lambda\cC)_p$.
\end{proof}

\begin{defn}[{\cite[Definition~4.5(i)]{SteinerEmbedding}}]
\label{Steineratomicdef}
Let $E$ be the basis of an $\omega$\nobreakdash-categ\-ory $\cC$. The basis $E$ is said to be \emph{atomic} if, for every $q > 0$, every $q$-cell $y$ of~$E$ and every $0 \leq p < q$,
we have
\[\supp([s_p(y)]_p)\cap\supp([t_p(y)]_p)=\emptyset.\]
    \end{defn}

\begin{prop}[{\cite[Proposition~4.6]{SteinerEmbedding}}]
 \label{remarkunital}  
 If $\cC$ is an $\omega$-category with
    atomic basis $E$, then the basis of the augmented directed chain complex $\lambda \cC$
    associated to $\cC$ is unital.
\end{prop}

\begin{defn}[{\cite[Definition~4.5(ii)]{SteinerEmbedding}}]
\label{defn:Steiner_loopfree_wcats}
Let $E$ be the basis of an $\omega$\nobreakdash-category $\cC$. The basis $E$ is said to be \emph{strongly loop-free} if the basis of the augmented directed chain complex $\lambda\cC$
is strongly loop-free in the sense of~\cref{defbasisalgebraic}.
\end{defn}

 \begin{defn}
\label{defnstrongSteiner}
An $\omega$-category $\cC$ is a \emph{strong Steiner $\omega$-category} if it admits an atomic and strongly loop-free basis.\footnote{Maltsiniotis and the first author in \cite[\textsection2.15]{AraMaltsiniotisJoin} define a \emph{strong Steiner $\omega$-category} to be any $\omega$-category in the essential image of strong Steiner complexes under $\nu$. In virtue of~\cite[Theorem~5.11]{SteinerEmbedding}, this is equivalent to being a strong Steiner $\omega$-category in the sense of~\cref{defnstrongSteiner}.} 
\end{defn}

More generally, all $\omega$-categories discussed in~\cref{ex:disks,ex:sphere,ex:theta,ex:oriental} are strong Steiner $\omega$-categories.
 For more details we refer the reader to \cite[Example~3.8]{SteinerEmbedding} for orientals, and to \cite{SteinerSimpleOmega}
 for the case of elements of $\Theta_n$. 
 
 \medskip

The reader shall notice that the conditions that define strong Steiner $\omega$\nobreakdash-categories make use of the associated augmented directed chain complex functor $\lambda$ from \cref{sec:nu}. The purpose of our paper is to give an intrinsic characterization of strong Steiner $\omega$-categories, proved as \cref{equivalencestrongSteiner}, that in particular does not rely on the functor $\lambda$.

\subsection{Steiner's Theorem}
\label{SteinersThm}
We can now recall the core result of Steiner's theory, which asserts a correspondence between strong Steiner $\omega$-categories and strong Steiner complexes.

\begin{thm}[{\cite[\textsection5]{SteinerEmbedding}}]
\label{nuFullyFaithful}
Let $\cC,\cC'$ be $\omega$-categories, and $C,C'$ augmented directed chain complexes.
\begin{enumerate}[leftmargin=*] 
    \item The complex $C$ is a strong Steiner complex if and only if $\nu C$ is a strong Steiner $\omega$\nobreakdash-category.
    \item  The $\omega$-category $\cC$ is a strong Steiner $\omega$-category if and only if $\lambda\cC$ is a strong Steiner complex.
    \item If $C$ and $C'$ are strong Steiner complexes, then $\nu$ induces a natural bijection
\[\adch(C,C')\cong \omega\cat(\nu C,\nu C').\]
    \item If $\cC$ and $\cC'$ are strong Steiner $\omega$-categories, then $\lambda$ induces a natural bijection
\[\omega\cat(\cC,\cC')\cong \adch(\lambda \cC,\lambda\cC').\]
\end{enumerate}
 \end{thm}

This means that the adjunction $\lambda\colon \omega\cat\rightleftarrows\adch\cocolon\nu$
restricts to an equivalence of categories
\begin{equation}
\label{SteinerEquivalence}
\lambda\colon\{\cC\text{ strong Steiner $\omega$-category}\}\simeq\{C\text{ strong Steiner complex}\}\cocolon\nu
\end{equation}
between the full subcategory of strong Steiner $\omega$-categories and the full subcategory of strong Steiner complexes. Given a strong Steiner $\omega$-category $\cC$, it is then justified to regard the augmented directed chain complex $C=\lambda\cC$ as an algebraic model of $\cC$.

\section{A categorical take on strong Steiner \pdfomega-categories}
\label{sec:Categorical}

The goal of this section is to characterize strong Steiner $\omega$-categories using categorical --- rather than algebraic --- terms, culminating in~\cref{equivalencestrongSteiner}.

\subsection{A categorical take on bases of \pdfomega-categories}

\label{sec:strong_Steiner}

We now discuss an established notion of freeness for an $\omega$-category. It was introduced independently by Street \cite[\S 4]{StreetOrientedSimplexes} under the name of computad, and by Burroni under the name of polygraph \cite[\textsection1.3]{BurroniWordProblems}.

Recall that we denote by $\varpi_{n}\colon\omegacat\to m\cat$ the truncation functor from \cref{truncation}, by $\mathcal D[n]$ the $n$-disk from~\cref{ex:disks}, and by
$\mathcal S[n-1]$ the $(n-1)$\nobreakdash-sphere $\varpi_{n-1} \mathcal D[n]$ from~\cref{ex:sphere}. Given an $\omega$-category $\cC$ and a set of cells $E$,
    for any $n>0$ one can consider a commutative square of the following form.
    \begin{equation}
    \label{attachcells}
        \begin{tikzcd}
            \coprod\limits_{E \cap \cC_n} \mathcal S[n-1] \ar[r] \ar[d,hook]    & \varpi_{n-1} \cC \ar[d,hook] \\
            \coprod\limits_{E \cap \cC_n} \mathcal D[n] \ar[r]           & \varpi_n\cC
        \end{tikzcd}
            \end{equation}
  Here, the vertical maps are the canonical inclusions, the horizontal component $\cD[n]\to\varpi_n\cC$ corresponding to $f\in E\cap\cC_n$ is $f$, and the horizontal component $\cS[n-1]\to\varpi_{n-1}\cC$ corresponding to $f\in E\cap\cC_n$ is the $(n-1)$\nobreakdash-dimensional boundary of $f$.

\begin{defn}
\label{defnpolygraph}
    Let $\cC$ be an $\omega$-category and $E$ a set of cells of $\cC$.
    We say that $E$ is a \emph{polygraphic generating set} of $\cC$, or that $\cC$ is \emph{freely generated}
    by~$E$, if 
    we have that $\cC_0\subseteq E$ and for any $n>0$ the commutative square~(\ref{attachcells})
    expresses the $\omega$-category $\varpi_n\cC$ as a pushout of $\omega$-categories.
\end{defn}

 Roughly speaking, the condition that $\cC$ be freely generated by a polygraphic generating set $E$ requests that the $n$-dimensional cells of $\cC$ are freely generated by the cells in $E\cap\cC_n$, which are attached along their boundary on~$\varpi_{n-1}\cC$.

\begin{prop}[{\cite[Proposition 1.5]{AraMaltsiniotisJoin}}]
\label{PolygraphImpliesCompositionGenerates}
    If an $\omega$-category $\cC$ is freely generated by a polygraphic generating set $E$, then $E$ composition-generates~$\cC$.
\end{prop}

  While it may seem like there could be different choices of polygraphic generating sets~$E$ that freely generate $\cC$, the following proposition identifies an intrinsic property in terms of the categorical structure of $\cC$ that characterizes the cells of $\cC$ that belong to $E$. This discussion is treated e.g.~in  \cite[\textsection5]{MakkaiNotes} and \cite[\textsection1.2.3]{HenryNonunitalPolygraphs}.
  
  \begin{prop}
  \label{generatorsareindecomposable}
  If an $\omega$-category is freely generated by a polygraphic generating set~$E$, the set $E$ is given precisely by all non-trivial cells of $\cC$ that are indecomposable, i.e., that can only be factored by using trivial cells. 
  \end{prop}

  In particular, if $\cC$ is freely generated by a polygraphic generating set $E$, then the set~$E$ is uniquely determined, and it makes sense to talk about \emph{the} polygraphic generating set that generates $\cC$. Notice that the uniqueness of the polygraphic generating set will also be a consequence of~\cref{polygraphimpliesbasis} together with~\cref{UniqueOmegaBasis}.
  
  All examples from~\cref{ex:disks,ex:sphere,ex:oriental,ex:theta} are freely generated, with the evident polygraphic generating sets.

One may regard~\cref{defnpolygraph} as a freeness condition that competes with Steiner's requirement of the existence of a basis (see~\cref{defnbasis}). While the two viewpoints are undoubtedly related, the two notions are not equivalent. On the one hand, one implication always holds.

\begin{prop}
\label{polygraphimpliesbasis}
If an $\omega$-category $\cC$ is freely generated by a polygraphic generating set
$E$, then $E$ is a basis of $\cC$.
\end{prop}

The main ingredients to prove the proposition are already treated as \cite[Proposition~4.2.4]{HenryRegular}.

\begin{proof}
To start with, the set of cells $E$ composition-generates $\cC$ by \cref{PolygraphImpliesCompositionGenerates}.
Now
the $\omega$-category $\cC$ is the colimit of the tower
    \[
     \varpi_0 \cC \to \varpi_1 \cC \to \cdots \to \varpi_q \cC \to \cdots 
    \]
    so, by applying
    the left adjoint functor $\lambda$ (from~\cref{sec:nu}), we obtain that
    $\lambda\cC$ is the colimit of the tower
     \[
     \lambda\varpi_0 \cC \to \lambda\varpi_1 \cC \to \cdots \to \lambda\varpi_q \cC \to \cdots 
    \]

    We consider the truncation operator $\varpi_q\colon\adch\to\adch$, which assigns to each augmented directed chain complex a new one that is trivial in degree higher than $q$. This operator $\varpi_q$ can be found in \cite[Chapter~2]{AraMaltsiniotisJoin} under the notation $\tau_{\leq q}^{\mathrm{b}}$. By \cite[Proposition~2.22]{AraMaltsiniotisJoin}, this chain complex truncation operator $\varpi_q$ commutes with the linearization functor $\lambda$, namely $\lambda\varpi_q \cC \cong \varpi_q\lambda \cC$. Then, since $\varpi_q\cC$ is characterized by the pushout (\ref{attachcells}), we obtain that the following pushout of augmented directed chain complexes.

\[
 \begin{tikzcd}
            \bigoplus\limits_{E \cap \cC_q} \lambda\mathcal S[q-1] \ar[r] \ar[d,hook]    & \lambda\varpi_{q-1} \cC \ar[d,hook]\ar[r,"\cong", phantom]& \varpi_{q-1}\lambda \cC\ar[d,hook]\\
            \bigoplus\limits_{E \cap \cC_q} \lambda\mathcal D[q] \ar[r]           & \lambda\varpi_q\cC\ar[r,"\cong", phantom]&\varpi_q\lambda\cC
        \end{tikzcd}
        \]
 By taking the group of $q$-chains, we obtain the pushout of abelian groups.
\begin{center}
\begin{tikzpicture}[scale=2]
   \tikzstyle{column 1}=[anchor=base]
  \tikzstyle{column 2}=[anchor=base]
    \matrix[nodes={},
        row sep={2cm,between origins},column sep={5cm,between origins},] {
          \node[strike out, draw](a1){ $\bigoplus\limits_{E \cap \cC_q} (\lambda\mathcal S[q-1])_q$ }; &
          \node[strike out, draw](a2){$(\lambda\varpi_{q-1} \cC)_q \cong
          (\varpi_{q-1}\lambda \cC)_q$};\\
           \node(b1) { $\bigoplus\limits_{E \cap \cC_q} (\lambda\mathcal D[q])_q$}; & 
          \node(b2){$(\lambda\varpi_q\cC)_q 
          \cong
          (\varpi_q\lambda\cC)_q$};\\
          };
          \draw[->, yshift=0.1cm] (a1.base east)--(a2.base west) node[midway, above]{$\cong$};
          \draw[->, yshift=0.1cm] (b1.base east)--(b2.base west) node[midway, above]{$\cong$};
          \draw[right hook->] (a1.south)--(b1.north);
          \draw[right hook->] (a2.south)--(b2.north);
\end{tikzpicture}
\end{center}
        
In particular, we obtain the isomorphism of abelian groups:
\[(\lambda\cC)_q \cong(\varpi_q\lambda\cC)_q\cong\bigoplus_{x\in E\cap\cC_q}(\lambda\cD[q])_q \cong \bigoplus_{x\in E\cap\cC_q}\mathbb Z[x]
                \cong \mathbb Z[E\cap\cC_q].\]
       Moreover, by inspection we see that this isomorphism restricts to a bijection
       \[(\lambda\cC)_q^+\cong\mathbb N[E\cap\cC_q],\]
       and the assignment $e\mapsto[e]_q$ defines a bijection between $E\cap\cC_q$ and the set $[E\cap\cC_q] \coloneqq \{[e]_q\in(\lambda\cC)_q\mid e\in E\cap\cC_q\}$.
\end{proof}

On the other hand, the following example will be used to show that, perhaps counter-intuitively, not every $\omega$-category with basis is freely generated by a polygraphic generating set.
 
\begin{ex}
\label{exampleforest}
Let $\cA$ be the $3$-category considered by Forest in \cite[\S 1.4]{ForestUnifying}, which is generated by the data
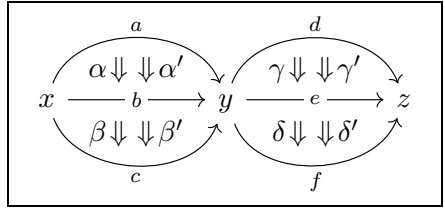
\[
\fbox{
  \begin{tikzcd}[ampersand replacement=\&]
    x \arrow[rr,"b"{description},""{auto=false,name=b}]
    \arrow[rr,out=70,in=110,"a",""{auto=false,name=a}] 
    \arrow[rr,out=-70,in=-110,"c"',""{auto=false,name=p}] 
    \arrow[phantom,"\alpha\!\Downarrow\ \Downarrow\!\alpha'",from=a,to=b]
    \arrow[phantom,"\beta\!\Downarrow\ \Downarrow\!\beta'",from=b,to=p] 
    \& \& y \arrow[rr,"e"{description},""{auto=false,name=e}]
    \arrow[rr,out=70,in=110,"d",""{auto=false,name=p}]
    \arrow[rr,out=-70,in=-110,"f"',""{auto=false,name=f}] \& \& z
    \arrow[phantom,"\gamma\!\Downarrow\ \Downarrow\!\gamma'",from=p,to=e]
    \arrow[phantom,"\delta\!\Downarrow\ \Downarrow\!\delta'",from=e,to=f]
  \end{tikzcd}
}
\]
together with the following $3$-cells.
\[
\begin{tikzcd}[ampersand replacement=\&]
    x \arrow[rr,"b"{description},""{auto=false,name=b}]
    \arrow[rr,out=70,in=110,"a",""{auto=false,name=a}] 
    \arrow[phantom,"\Downarrow \alpha",from=a,to=b]
    \& \& y \arrow[rr,"e"{description},""{auto=false,name=e}]
    \arrow[rr,out=-70,in=-110,"f"',""{auto=false,name=f}] \& \& z
    \arrow[phantom,"\Downarrow \delta",from=e,to=f]
  \end{tikzcd}
  \quad\overset{A}{\Rrightarrow}\quad
  \begin{tikzcd}[ampersand replacement=\&]
    x \arrow[rr,"b"{description},""{auto=false,name=b}]
    \arrow[rr,out=70,in=110,"a",""{auto=false,name=a}] 
    \arrow[phantom,"\Downarrow \alpha'",from=a,to=b]
    \& \& y \arrow[rr,"e"{description},""{auto=false,name=e}]
    \arrow[rr,out=-70,in=-110,"f"',""{auto=false,name=f}] \& \& z
    \arrow[phantom,"\Downarrow \delta'",from=e,to=f]
  \end{tikzcd}
  \]
  \[
    \begin{tikzcd}[ampersand replacement=\&]
    x \arrow[rr,"b"{description},""{auto=false,name=b}]
    \arrow[rr,out=-70,in=-110,"c"',""{auto=false,name=p}]
    \arrow[phantom,"\Downarrow \beta",from=b,to=p]
    \& \& y \arrow[rr,"e"{description},""{auto=false,name=e}]
    \arrow[rr,out=70,in=110,"d",""{auto=false,name=q}]
    \& \& z
    \arrow[phantom,"\Downarrow \gamma",from=q,to=e]
  \end{tikzcd}
  \quad\overset{B}{\Rrightarrow}\quad
  \begin{tikzcd}[ampersand replacement=\&]
    x \arrow[rr,"b"{description},""{auto=false,name=b}]
    \arrow[rr,out=-70,in=-110,"c"{description},""{auto=false,name=p}] 
    \arrow[phantom,"\Downarrow \beta'",from=b,to=p]
    \& \& y \arrow[rr,"e"{description},""{auto=false,name=e}]
    \arrow[rr,out=70,in=110,"d",""{auto=false,name=q}]
    \& \& z
    \arrow[phantom,"\Downarrow \gamma'",from=q,to=e]
  \end{tikzcd}
  \]
  
It is shown in loc.~cit.~that the $3$-category $\cA$ has 
two distinct parallel $3$-cells $H_1$ and $H_2$ that are defined as different compositions involving the
    same $3$-cells of the basis: \[H_1\coloneqq\big((a*_0\gamma)*_1A*_1(\beta*_0f)\big)*_2\big((\alpha'*_0d)*_1B*_1(c*_0\delta')\big)\]
    and
    \[H_2\coloneqq\big((\alpha*_0d)*_1B*_1(c*_0\delta)\big)*_2\big((a*_0\gamma')*_1A*_1(\beta'*_0f)\big).\]
Notice that we have $[H_1]_3 = [H_2]_3 = [A]_3 + [B]_3$ in $(\lambda \cA)_3$.
Let $\cB$ be the $3$-category obtained by identifying the two cells $H_1$ and $H_2$, namely defined by the following pushout, in which the left vertical map is the folding map of $\cD[3]$.
\[
\begin{tikzcd}
            \cD[3]\amalg\cD[3] \ar[rr,"{[H_1,H_2]}"] \ar[d, two heads]    && \cA \ar[d, two heads] \\
            \cD[3] \ar[rr]           && \cB
        \end{tikzcd}
\]
It turns out that $\cB$ is a counter-example to the converse implication of~\cref{polygraphimpliesbasis}.
\end{ex}

\begin{prop}
The $3$-category $\cB$ from~\cref{exampleforest} has the following properties:
\begin{enumerate}[leftmargin=*]
    \item $\cB$ admits the basis
    \[E\coloneqq\{x,y,z,a,b,c,d,e,f,\alpha, \alpha', \beta, \beta', \gamma, \gamma', \delta,\delta', A, B\};\]
    \item $\cB$ is not freely generated by a polygraphic generating set.
\end{enumerate}
\end{prop}

We shall give a quick sketch
    of this fact, since providing a full proof is beyond the scope of this work.

\begin{proof}[Sketch of the proof]
We start by proving (1). The 3-category $\cA$ is by definition freely generated by a polygraphic generating set, and by~\cref{polygraphimpliesbasis} we know that $\cA$ admits a basis.
Moreover, the canonical map $\cA\to\cB$ induces isomorphisms
\[(\lambda\cA)_q\cong(\lambda\cB)_q\cong\left\{
\begin{array}{cc}
    \mathbb Z[x,y,z] & q=0, \\
    \mathbb Z[a,b,c,d,e,f]   & q=1,\\
    \mathbb Z[\alpha,\alpha',\beta,\beta',\gamma,\gamma',\delta,\delta']   & q=2,\\
    0&q\ge 4.
\end{array}\right.\]
Since the classes $[H_1]_3$ and $[H_2]_3$ are equal in $(\lambda\cA)_3$
and the left adjoint functor $\lambda$ preserves pushouts, the
canonical map $\cA\to\cB$ also induces an isomorphism \[(\lambda\cA)_3\cong(\lambda\cB)_3\cong\mathbb Z[A,B].\]
In particular, the canonical map $\cA\to\cB$ induces isomorphisms $(\lambda\cA)_q\cong(\lambda\cB)_q$ for any $q\ge0$, and so $\cB$ admits a basis since $\cA$ does.

For (2), we observe that, if $\cB$ were freely generated by a polygraphic generating set, then, since  $\cA$ and $\cB$ have the same indecomposable cells, they both would be freely generated by the same set of cells
\[E=\{x,y,z,a,b,c,d,e,f,\alpha, \alpha', \beta, \beta', \gamma, \gamma', \delta,\delta', A, B\}.\]
In particular, we would have an isomorphism $\cA\cong\cB$ that preserves $H_1$ and~$H_2$. However, in $\cB$ we have $H_1=H_2$ by construction, while Forest shows in \cite[\S 1.4]{ForestUnifying} that in $\cA$ we have $H_1\neq H_2$.
\end{proof}

  \subsection{A categorical take on atomicity}
  
In this subsection we collect a series of considerations related to the atomicity condition (from~\cref{Steineratomicdef}) for a basis of an $\omega$-category freely generated by a polygraphic generating set.

The following definition of \emph{factor} leads to the definition of \emph{support} of a cell, which is inspired by Makkai's work \cite{MakkaiNotes}.

\begin{defn}
\label{def:factor_cell}
Let $\cC$ be an $\omega$-category freely generated by a polygraphic generating set~$E$.
We say that $x\in E$ is a \emph{factor}
of a $q$-cell $y$ if $y$ can be expressed as a composite
\[
y=a_1*_{k_1}*a_2*_{k_2}*\dots*a_r*_{k_r}x*_{k_{r+1}}a_{r+1}*_{k_{r+2}}\dots *_{k_{s-1}}a_{s-1}*_{k_{s}}a_{s}
\]
where $a_i$ is a cell in $\cC$ and $k_i\ge0$ for $i=1,\dots,s$, for some bracketing of the expression so that it makes sense in $\cC$.
\end{defn}

\begin{defn}
\label{context}
   Let $\cC$ be an $\omega$-category freely generated by a polygraphic generating set $E$. Given a $q$-cell $y$, the \emph{support} of $y$ is the (possibly empty\footnote{The support of a cell is empty if and only if the cell is trivial.}) set
   \[\supp(y) \coloneqq \{e\in E\cap\cC_{q}\mid \text{$e$ is a factor of $y$}\}\subseteq E\cap\cC_q.\]
\end{defn}

The following lemma explains how this notion of support is compatible with the one given in the algebraic context from~\cref{SupportComplex}.

\begin{lem}
   \label{rmk:equivalence_support}
   Let $\cC$ be an $\omega$-category freely generated by a polygraphic generating set\footnote{In fact, the given proof only uses that $E$ is a basis for the $\omega$-category $\cC$.}~$E$. Then, for any $q$-cell $y$ in $\cC$, the assignment
   \[ e\in E\cap\cC_q\mapsto[e]_q\in [E\cap\cC_q] \]
   induces a bijection
   \[\supp(y)\cong\supp([y]_q).\]
\end{lem}

\begin{proof}
For $q\ge0$, by~\cref{polygraphimpliesbasis}
the assignment $e\mapsto[e]_q$ is a bijection from $E\cap \cC_q$ to
the basis $[E \cap \cC_q]$
of the free abelian group $(\lambda \cC)_q$.

Let $y$ be a $q$-cell of $\cC$. By \cref{PolygraphImpliesCompositionGenerates}, $E$ composition-generates $\cC$, so there exists an expression $D$ of $y$ in terms of elements of $E$, which generally involves cells in dimension $q$ or lower, and we fix one such specific decomposition. Let $E_y(D)$ denote the elements of $E\cap\cC_q$ that occur in the fixed decomposition of $y$, and for any $e$ in $E_y(D)$ denote by $\lambda_e\geq 1$ the number of times $e$ appears
in the expression $D$. According to these definitions, the set $E_y(D)$ is
always a subset of $\supp(y)$.

We now show that the assignment $e\mapsto[e]_q$ defines a bijection
\[
E_y(D)\cong\supp([y]_q).
\]
By definition of $(\lambda\cC)_q$ and using~\cref{propertieslambda}, we obtain that the unique expression of $[y]_q$ of the form \eqref{canonicaldecomposition} is given by
\begin{equation}
\label{firstexpressiony}
    [y]_q = \sum_{e \in E_y(D)} \lambda_e \cdot [e]_q.
\end{equation}
Since $E_y(D)\subseteq\supp(y)\subseteq E\cap\cC_q$, every $[e]_q$ is an element of the basis of $(\lambda\cC)_q$,
so this expression is the unique decomposition of
$[y]_q$ as a linear combination of
basis elements in $(\lambda\cC)_q$ from \eqref{canonicaldecomposition}. In particular, the set $E_y(D)$ is canonically in bijection with the set $\supp([y]_q)$ via the assignment $e\mapsto[e]_q$.

This shows in particular that $E_y(D)$ does not depend on $D$.
Indeed, if $D'$ is another expression of $y$ in terms of elements of $E$,
the map $e \mapsto [e]_q$ from $E \cap \cC_q$ to $[E \cap \cC_q]$ restricts to a bijection $E_y(D') \cong \supp([y]_q)$, from which necessarily $E_y \coloneqq E_y(D) = E_y(D')$.
Hence, we obtain the converse inclusion $\supp(y) \subseteq E_y$.

In conclusion, we showed $\supp(y) = E_y \cong \supp([y]_q)$ via $e\mapsto[e]_q$, as desired.
\end{proof}

\begin{rmk}
\label{rmk:inclusion_support}
Let $\cC$ be an $\omega$-category freely generated by a polygraphic generating set~$E$. For any $q$-cell $y$ in $\cC$, all the coefficients
    occurring in the decompositions \eqref{canonicaldecomposition} of $[s_{q-1}(y)]_{q-1}$ and $[t_{q-1}(y)]_{q-1}$
    in $(\lambda \cC)_{q-1}$ are non-negative by \cref{positivecoefficients}. Moreover, by definition we have \[\partial([y]_q)=[t_{q-1}(y)]_{q-1}-[s_{q-1}(y)]_{q-1}.\]
    In particular, it follows that there are inclusions
   \[\supp (\partial^-([y]_q))\subseteq\supp([s_{q-1}(y)]_{q-1}) \]
   and
   \[ \supp (\partial^+([y]_q))\subseteq\supp([t_{q-1}(y)]_{q-1}).\]
\end{rmk}

The following example, closely related to~\cite[Exemple~3.4]{AraMaltsiniotisJoin}, shows that the inclusion of supports from~\cref{rmk:inclusion_support} is generally strict.

\begin{ex}
\label{ex:different_loop-freeness}
 Let $\cC$ be the following $2$-category.
 \[\cC=\fbox{\begin{tikzcd}[ampersand replacement=\&]
            x \ar[r, "f"] \ar[d, "f"']  \& y \ar[d, "g"] \\
            y \ar[r, "h"']              \& z
            \ar[Rightarrow, from=1-2, to=2-1, shorten <=3mm, shorten >=3mm, "\alpha"]
        \end{tikzcd}}\]
        Then, one can check that $\cC$ has the following properties:
\begin{enumerate}[leftmargin=*]
    \item The $\omega$-category $\cC$ is freely generated by the polygraphic generating set
    \[ E=\{x,y,z,f,g,h,\alpha\}; \]%
    in particular by~\cref{polygraphimpliesbasis}, we know that $E$ is a basis for $\cC$.
    \item We have 
\[\supp([s_{1}(\alpha)]_{1})\cap\supp([t_{1}(\alpha)]_{1})=\{[f]_1\}\neq\emptyset;\]
    in particular, $E$ is a non-atomic basis for $\cC$.
    \item There is a strict inclusion of supports
    \[\supp (\partial^+([\alpha]_2))=\{[h]_1\}\subsetneq\{[f]_1,[h]_1\}=\supp([t_{1}(\alpha)]_1).\]
\end{enumerate}
\end{ex}

\begin{lem}
   \label{lem:equivalence_support3}
   Let $\cC$ be an $\omega$-category freely generated by a polygraphic generating set~$E$
   and $y$ a $q$-cell of $\cC$ for $q>0$. If
   \[\supp([s_{q-1}(y)]_{q-1})\cap\supp([t_{q-1}(y)]_{q-1})=\emptyset,\]
    then we have the equalities
    \[
        [s_{q-1}(y)]_{q-1} = \partial^-([y]_{q})
        \ \text{ and }\ 
        [t_{q-1}(y)]_{q-1} = \partial^+([y]_q)
    \]
 and in particular the supports are equal:
 \[\supp (\partial^-([y]_q))=\supp([s_{q-1}(y)]_{q-1})\]
 and
 \[ \supp (\partial^+([y]_q))=\supp([t_{q-1}(y)]_{q-1}).\]
\end{lem}

\begin{proof}
    By definition, we have
    \[\partial([y]_q)=[t_{q-1}(y)]_{q-1}-[s_{q-1}(y)]_{q-1},\]
    and by~\cref{positivecoefficients} the coefficients appearing in the linear decompositions
 of $[t_{q-1}(y)]_{q-1}$ and $[s_{q-1}(y)]_{q-1}$ are all non-negative. The assumption that the supports of $[t_{q-1}(y)]_{q-1}$ and $[s_{q-1}(y)]_{q-1}$ are disjoint prevents the possibility that there could be common chains occurring in the linear decompositions
 of $[t_{q-1}(y)]_{q-1}$ and $[s_{q-1}(y)]_{q-1}$.
As a consequence, we obtain the equalities $[s_{q-1}(y)]_{q-1} = \partial^-([y]_{q})$
  and $[t_{q-1}(y)]_{q-1} = \partial^+([y]_q)$.
\end{proof}

\begin{rmk}
In particular, the lemma applies to any element of a polygraphic generating set which is atomic.
\end{rmk}

We now define two preorder relations $\preceq_{\mathbb N}^\omega$ and $\leq_{\mathbb N}^\omega$ on the polygraphic generating set of an
$\omega$-category. The relation $\preceq_{\mathbb N}^\omega$ is the naive generalization of Steiner's relation $\leq_{\mathbb N}$ from
\cref{loopfree_preorder_complexes}. It will turn out to be insufficient to capture the correct $\omega$-categorical analog of Steiner's notion of a strong loop-free basis. This is the reason why we introduce the relation $\leq_{\mathbb N}^\omega$.
We'll see in \cref{lem:equivalence_w-relations} that in presence of atomicity the two relations agree.

\begin{notn}
\label{preorder_omegacat_codim1}
Let $\cC$ be an $\omega$-category freely generated by a polygraphic generating set~$E$. We denote by $\preceq_{\mathbb N}^\omega$ the preorder relation on $E$ generated by the condition
\[\text{$a\preceq^\omega_{\mathbb N}b$ if $a\in E\cap\cC_p$, $b\in E\cap\cC_q$, and
$\begin{cases}
\text{either $a\in \supp(s_{q-1}(b))$ with $q > 0$,}\\
\text{or $b\in \supp (t_{p-1}(a))$ with $p > 0$.}
\end{cases}$}
\]
\end{notn}

\begin{notn}
\label{preorder_omegacat}
Let $\cC$ be an $\omega$-category freely generated by a polygraphic generating set $E$. We denote by $\leq_{\mathbb N}^\omega$ the preorder relation on $E$ generated by the condition
\[\text{$a\leq_{\mathbb N}^\omega b$ if $a \in E \cap \cC_p$, $b \in E \cap\cC_q$, and 
$\begin{cases}
\text{either $a\in \supp(s_{p}(b))$,}\\
\text{or $b\in \supp (t_{q}(a))$.}
\end{cases}$}
\]
\end{notn}

\begin{rmk}
\label{lem:equivalence_w-relations_easy}
    Let $\cC$ be an $\omega$-category freely generated by a polygraphic generating set~$E$.
    Then the preorder $\leq_\mathbb{N}^\omega$ is finer than the preorder $\preceq^\omega_\mathbb{N}$
    on $E$. Namely, if $a \preceq^\omega_\mathbb{N} b$, then $a \leq_\mathbb{N}^\omega b$ for $a,b\in E$.
\end{rmk}

After a preliminary fact, we show in~\cref{lem:equivalence_w-relations} that in presence of atomicity the two relations $\preceq_\mathbb N^\omega$ and $\leq_\mathbb N^\omega$ in fact coincide.

\begin{lem}
\label{precandSteiner}
Let $\cC$ be an $\omega$-category freely generated by a polygraphic generating set~$E$. 
\begin{enumerate}[leftmargin=*]
    \item The assignment
    $e\in E\mapsto[e]\in[E]:=\coprod_{q\ge0}[E\cap\cC_q]$
    is preorder-creating for the preorders
    $\preceq^\omega_\mathbb{N}$ in the source and $\leq_{\mathbb{N}}$ in the target. 
    Namely, if $[a]\leq_{\mathbb N}[b]$ then $a\preceq_{\mathbb N}^\omega b$ for $a,b\in E$. 
    \item If $E$ is an atomic basis, then the bijection is in addition order-preserving.
    Namely, in presence of atomicity, $[a]\leq_{\mathbb N}[b]$ if and only if $a\preceq_{\mathbb N}^\omega b$ for~$a,b\in E$.
\end{enumerate}
\end{lem}

\begin{proof}
For the first part we observe that, if $a\in E\cap\cC_p$, $b\in E\cap\cC_q$ and $[a]_p \in\supp(\partial^-([b]_q))$,
then $q=p+1$ and by
\cref{rmk:inclusion_support} we obtain
\[[a]_p\in\supp(\partial^-([b]_{p+1}))\subseteq\supp([s_p(b)]_p).\]
Then by~\cref{rmk:equivalence_support} we get
\[a\in\supp(s_p(b)).\]
In particular, it follows that $a\preceq_{\mathbb N}^\omega b$. A similar argument applies to the case $[b]_q \in\supp(\partial^+([a]_p))$.

For the second part we observe that, if $a\in E\cap\cC_p$, $b\in E\cap\cC_q$ and $a\in\supp(s_p(b))$, then $q=p+1$ and by atomicity combined with~\cref{lem:equivalence_support3} we obtain
\[[a]_p\in\supp([s_p(b)]_p)=\supp(\partial^-([b]_{p+1})).\]
In particular, it follows that $[a]_p\leq^{\omega}_\mathbb N[b]_{p+1}$.
A similar argument applies to the case $b \in\supp(t_p(a))$.
\end{proof}

 \begin{prop}
\label{lem:equivalence_w-relations}
Let $\cC$ be an $\omega$-category freely generated by a 
polygraphic generating set~$E$. 
   If $E$ is an atomic basis, then the preorders $\preceq^\omega_\mathbb{N}$
    and $\leq_\mathbb{N}^\omega$ on~$E$ coincide. Namely, in presence of atomicity $a\leq_\mathbb{N}^\omega b$ if and only if $a\preceq_\mathbb{N}^\omega b$ for $a,b\in E$.
\end{prop}

\begin{proof}
We know by~\cref{lem:equivalence_w-relations_easy} that,
    if $a \preceq^\omega_\mathbb{N} b$, then $a \leq_\mathbb{N}^\omega b$, and we now prove the converse under the assumption that $E$ is atomic.

  For this, we suppose that $a\in E\cap\cC_p$, $b\in E\cap\cC_q$ and $a \in\supp(s_p(b))$, which imposes the constraint $q-p\ge1$, and we prove by induction on $q-p\ge1$ that $a \preceq^\omega_\mathbb{N} b$. The other case, namely $b \in\supp(t_q(a))$, can be treated analogously. Overall, this will have shown that if $a\leq_\mathbb N^\omega b$, then $a\preceq_\mathbb N^\omega b$.
  
If $q-p=1$, then by definition we have
 \[a\in\supp(s_{p}(b))=\supp(s_{q-1}(b))\]
 and in particular \[a\preceq_\mathbb N^\omega b.\]
If $q-p>1$, the atomicity condition
guarantees that
\[\supp(s_ps_{p+1}(b))\cap\supp(t_ps_{p+1}(b))=\emptyset,\]
so we can apply~\cref{lem:equivalence_support3}
in the equalities below:
\[\begin{array}{cllll}
  [a]_p&\in&\supp([s_p(b)]_p) &\text{\cref{rmk:equivalence_support}}\\   
  &=&\supp([s_ps_{p+1}(b)]_p) &\text{\eqref{globularity}}\\
  &=& \supp(\partial^-[s_{p+1}b]_{p+1}) &\text{\cref{lem:equivalence_support3}} \\
&\subseteq&   \bigcup\limits_{e\in \supp [s_{p+1}(b)]_{p+1}}\supp(\partial^-([e]_{p+1})).  &
\end{array}\]
This says that there exists $c \in \supp(s_{p+1}(b))$ such that
\[[a]_p\in\supp(\partial^-([c]_{p+1})).\]
By~\cref{rmk:equivalence_support,lem:equivalence_support3},
this means that
\[a\in \supp(s_{p}(c)).\]
Given that
  \[a\in\supp(s_p(c))\quad\text{ and }\quad c\in\supp(s_{p+1}(b)),\]
   we have
    \[a\preceq_{\mathbb N}^\omega c\quad\text{ and }\quad c\leq_{\mathbb N}^\omega b.\]
  By induction hypothesis applied to the pair $c$ and $b$,
    this implies
    \[a\preceq_{\mathbb N}^\omega c\preceq_{\mathbb N}^\omega b,\]
as desired.
\end{proof}

As a consequence, we obtain that under the hypothesis of atomicity the relation $\leq_{\mathbb N}^\omega$ is the correct analog of Steiner's relation
$\leq_{\mathbb N}$:

 \begin{cor}  
 \label{lemmapreorder2}
 If an $\omega$-category $\cC$ is freely generated by a polygraphic generating set $E$ that is atomic in the sense of~\cref{Steineratomicdef}, then the assignment $e\in E\mapsto[e]\in[E]$
    is an isomorphism of pre-ordered sets from $(E, \leq_{\mathbb N})$ to~$([E], \leq_{\mathbb N}^\omega)$. Namely, under these circumstances, if $[a]\leq_{\mathbb N}[b]$ then $a\leq_{\mathbb N}^\omega b$ for~$a,b\in E$.
\end{cor}

\begin{proof}
This is a direct consequence of~\cref{precandSteiner,lem:equivalence_w-relations}.
\end{proof}

\subsection{A categorical take on strong loop-freeness}

We are now ready to propose a categorical notion of strong loop-freeness for $\omega$-categories freely generated by a polygraphic generating set, which we'll show in \cref{rmk:loop-free,prop:atomicity+loop-free} to be compatible with the algebraic one from \cref{defn:Steiner_loopfree_wcats}.

\begin{defn}
\label{stronglyloopfreepolygraph}
    Let $\cC$ be an $\omega$-category freely generated by a polygraphic generating set~$E$.
    The polygraphic generating set $E$ is said to be \emph{strongly loop-free} if 
    the preorder relation (namely, reflexive and transitive relation) $\leq_{\mathbb N}^\omega$ on $E$ defines in fact a partial order, i.e., if $\leq_{\mathbb N}^\omega$ is antisymmetric.
\end{defn}

 To give some familiarity with the relation $\leq_{\mathbb N}^\omega$ and with the notion of strong loop-freness, let's look at a couple of examples.
\begin{ex}
Consider the $2$-oriental $\cO[2]$,
\[
	\cO[2] = 
		\fbox{	 \begin{tikzcd}[row sep=small, column sep=tiny,ampersand replacement=\&]
                                            \& b \ar[rd,"g"]   \&         \\
       a \ar[ru,"f"] \ar[rr, "h"{name=s, below}] \&                   \& c
                \ar[Rightarrow, from=s, to=1-2, shorten >= 2pt, shorten <= 2pt, "\tiny{\alpha}"]
\end{tikzcd}
}
\]
which is freely generated by the (atomic) polygraphic generating set 
\[E=\{a,b,c,f,g,h,\alpha\}.\]
The relation $\leq_{\mathbb N}^\omega$ on $E$ is in fact a total order, given by
\[a\leq_{\mathbb N}^\omega  h\leq_{\mathbb N}^\omega \alpha\leq_{\mathbb N}^\omega f \leq_{\mathbb N}^\omega b\leq_{\mathbb N}^\omega g \leq_{\mathbb N}^\omega c,\]
and $\cO[2]$ is therefore freely generated by the strongly loop-free polygraphic generating set $E$.
\end{ex}

More generally, all $\omega$-categories discussed in~\cref{ex:disks,ex:sphere,ex:theta,ex:oriental} are freely generated by strongly loop-free and atomic polygraphic generating sets.

On the contrary, the following example
illustrates how the condition of strong loop-freeness prevents the existence of endomorphisms.

\begin{ex}
 \label{ex:loop}
Consider the $1$-category $\cC$,
\[
\cC=
\fbox{\begin{tikzcd}[ampersand replacement=\&]
				a \arrow[r, bend left=50, "f"{name=U, above}] 
				\& b\arrow[l, bend left=50, "g"{name=D}]
			\end{tikzcd}
			}
			\]
which is freely generated by a polygraphic generating set $E=\{a,b,f,g\}$. We have
\[a\leq_{\mathbb N}^\omega f\leq_{\mathbb N}^\omega b\leq_{\mathbb N}^\omega g\leq_{\mathbb N}^\omega a,\]
while $a\neq b$. Hence $\cC$ is not freely generated by a strongly loop-free
polygraphic generating set.
\end{ex}

The following (non-)example illustrates a different phenomenon that the strong loop-freeness condition is designed to prevent.

\begin{ex} 
Consider the $2$-category $\cA$ corresponding to the following picture.
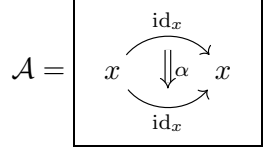
\[
			\cA = 
			\fbox{
			\begin{tikzcd}[ampersand replacement=\&]
				x \arrow[r, bend left=50, ""{name=U, below}, "\id_x"{above}] \arrow[r, bend right=50, ""{name=D}, "\id_x"{below}]
				\& x
				\ar[Rightarrow, from=U, to=D, shorten >= 2pt, "\alpha"]
			\end{tikzcd}
			}
			\]
			Equivalently, $\cA$ has a single $0$-cell $x$, no non-trivial $1$-cells and is generated by a single non-trivial $2$-cell $\alpha$.
In particular, $\cA$ is freely generated by the polygraphic generating set $\{x,\alpha\}$,
which is evidently \emph{not} atomic. The preorder relation $\preceq_{\mathbb N}^\omega$ on
$\{x, \alpha\}$ is discrete and it is therefore a partial order relation.
However, since
\[x\in \supp(s_0(\alpha))\cap\supp(t_0(\alpha)),\]
we have
\[\alpha\leq_\mathbb N^\omega x\leq_\mathbb N^\omega\alpha,\]
although $\alpha\neq x$. Hence $\cA$  is not freely generated by a strongly loop-free
polygraphic generating set. This shows that without atomicity of the basis
the relation $\preceq_{\mathbb N}^\omega$ can be a partial order, while $\leq_{\mathbb N}^\omega$
is not. In particular, these two preorder relations may differ in the absence of atomicity. 
\end{ex}

A crucial consequence of the notion introduced in \cref{stronglyloopfreepolygraph} is atomicity.

\begin{prop}
 \label{rmk:atomic}
    If an $\omega$-category $\cC$ is freely generated by a strongly loop-free polygraphic generating set $E$,
    then $E$ is an atomic basis of $\cC$.
    \end{prop}

\begin{proof}
     Let $e$ be an element of $E\cap \cC_q$, for $q >0$. We wish to prove that for $0 \le p <q$ \[\supp([s_p(e)]_p)\cap\supp([t_p(e)]_p)=\emptyset.\]
     By~\cref{rmk:equivalence_support},
     this is equivalent to proving that
     \[\supp(s_p(e))\cap\supp(t_p(e))=\emptyset.\]
Now, assuming there exists a cell
\[b\in\supp(s_p(e))\cap\supp(t_p(e)),\]
     we would have the relations \[b \leq_{\mathbb{N}}^{\omega} e \leq_{\mathbb{N}}^\omega b,\]
     so by strong  loop-freeness $b=e$, which is impossible because $e$ is a $q$-cell and $b$ is a $p$-cell and $p<q$.
 \end{proof}

We are now ready to prove that the algebraic and categorical notions of strong loop-freeness agree in presence of atomicity.

\begin{prop}
 \label{rmk:loop-free}
    If an $\omega$-category $\cC$ is freely generated by a strongly loop-free polygraphic generating set $E$,
    then $E$ is a strongly loop-free basis of $\cC$.
    \end{prop}

    \begin{proof}
    The fact that $E$ is a basis for $\cC$ follows from~\cref{polygraphimpliesbasis}, the fact that $E$ is atomic follows from \cref{rmk:atomic}, and the fact that the basis $E$ is strongly loop-free from~\cref{lemmapreorder2}.
   \end{proof}

It is not true that if $\cC$ is freely generated by a polygraphic generating set~$E$ and $E$ is a strongly loop-free basis for $\cC$ then $\cC$ is freely generated by a strongly loop-free polygraphic generating set $E$.
This can be argued using the following example, which is closely
    related to~\cite[Exemple~3.4]{AraMaltsiniotisJoin}. We leave the verifications of the details to the interested reader.
\begin{ex}
\label{ex:different_loop-freeness2}
 Let $\cC$ be the $2$-category from~\cref{ex:different_loop-freeness}.
 \[\cC=\fbox{\begin{tikzcd}[ampersand replacement=\&]
            x \ar[r, "f"] \ar[d, "f"']  \& y \ar[d, "g"] \\
            y \ar[r, "h"']              \& z
            \ar[Rightarrow, from=1-2, to=2-1, shorten <=3mm, shorten >=3mm, "\alpha"]
        \end{tikzcd}}\]
Then, we saw that $\cC$ is freely generated by a polygraphic generating set $E$ and one can check that the basis $[E]$ is strongly loop-free. Nevertheless, 
$\cC$~is not freely generated by a strongly loop-free polygraphic generating set. Indeed, we have that 
$f\leq_{\mathbb N}^\omega\alpha\leq_{\mathbb N}^\omega f$, while $f\neq\alpha$.
\end{ex}

\cref{ex:different_loop-freeness2} shows that for an $\omega$-category freely generated by a polygraphic generating set
the ``$\omega$-categorical'' and the ``algebraic'' notions of strong loop-freeness
are not equivalent.
As~\cref{prop:atomicity+loop-free} shows, the lack of the atomicity condition was in fact the only obstruction.

\begin{prop}
    \label{prop:atomicity+loop-free}
     Let $\cC$ be an $\omega$-category freely generated
     by a polygraphic generating set~$E$.
     If $E$ is an atomic and strongly loop-free basis for $\cC$,
     then $\cC$ is freely generated by a strongly loop-free polygraphic generating set $E$.
\end{prop}

\begin{proof}
We prove that the relation $\leq_{\mathbb N}^{\omega}$ is antisymmetric. Given a $a$ in $E \cap \cC_p$ and $b$ in $E \cap \cC_q$, if
\[a\leq_{\mathbb N}^\omega b\quad\text{ and }\quad b\leq_{\mathbb N}^\omega a,\]
then by~\cref{lemmapreorder2}
we have
\[[a]_p\leq_{\mathbb N}[b]_q\quad\text{ and }\quad[b]_q\leq_{\mathbb N}[a]_p.\]
Since the basis $[E]$ of $\lambda\cC$ is strongly loop-free, we deduce
\[[a]_p=[b]_q,\]
and by~\cref{defnbasis} we deduce
\[a=b,\]
which concludes the proof.
\end{proof}

\subsection{A categorical take on strong Steiner \pdfomega-categories}

The following theorem, finally, establishes how our notion of strongly loop-free polygraphic generating set is an alternative description for the same class of $\omega$-categories considered in Steiner's theory.

\begin{thm}
\label{equivalencestrongSteiner}
For an $\omega$-category $\cC$, the following are equivalent:
\begin{enumerate}[leftmargin=*,label=(\arabic*)]
    \item\label{omega} $\cC$ is freely generated by a strongly loop-free polygraphic generating set;
    \item\label{strongsteiner} $\cC$ is a strong Steiner $\omega$-category, i.e., it admits a strongly loop-free atomic basis;
    \item\label{nuvalue} $\cC$ is isomorphic as an $\omega$-category to $\nu C$ for some strong Steiner complex~$C$.
\end{enumerate}
\end{thm}

\begin{proof}
The equivalence between \ref{strongsteiner} and \ref{nuvalue} is Steiner's ~\cref{nuFullyFaithful}, and we now show the equivalence of \ref{omega} and \ref{strongsteiner}.

We first prove that \ref{omega} implies \ref{strongsteiner}. Let $\cC$ be an $\omega$-category that is freely generated by a strongly loop-free polygraphic generating set $E$.
By~\cref{polygraphimpliesbasis}, we know that $\cC$ admits a basis $E$, by~\cref{rmk:loop-free} that the basis~$E$ is strongly loop-free, and by~\cref{rmk:atomic} that it is atomic.

We now prove that \ref{strongsteiner} implies \ref{omega}.
Let $\cC$ be an $\omega$-category admitting a strongly loop-free atomic basis $E$. 
Hence, using~\cref{remarkunital}, we know that
the associated augmented directed chain complex $\lambda\cC$ admits a strongly loop-free unital
basis.
Hence, by virtue of~\cref{nuFullyFaithful} the unit $\eta_{\cC}$ of the adjunction $(\lambda,\nu)$ is an isomorphism of $\omega$-categories 
$\eta_{\cC}\colon\cC\cong \nu \lambda \cC$. Finally, by~\cite[Theorem~6.1]{SteinerEmbedding} the $\omega$-category $\cC\cong\nu\lambda\cC$ is freely generated by a polygraphic generating set. We then conclude using~\cref{prop:atomicity+loop-free}.
\end{proof}

\begin{rmk}
  One can prove in essentially the same way an analogous theorem for another
  class of $\omega$-categories playing an important role in Steiner's theory
  \cite{SteinerEmbedding}, named Steiner $\omega$-categories in
  \hbox{\cite[\textsection2]{AraMaltsiniotisJoin}}. More precisely, one can
  prove that an $\omega$\nobreakdash-category $\mathcal{C}$
  is a Steiner $\omega$\nobreakdash-category if and only if it is freely
  generated by a polygraphic generating set $E$ on which there exists an order relation satisfying the following property:
  for every $e$ in $E \cap \mathcal{C}_p$ and every $q$ such 
  that $0 \le q \le p$, then each element of the support of $[s_q(e)]_q$ is
  strictly smaller than each element of the support of $[t_q(e)]_q$.
\end{rmk}

\bibliographystyle{amsalpha}
\bibliography{ref} 

\end{document}

%% file: preamble.tex
\usepackage{amsmath}
\usepackage{latexsym,amsfonts,amssymb,mathrsfs}
\usepackage{mathdots}
\usepackage{color}
\usepackage[all,2cell,color]{xy}
\usepackage{enumitem} 
\usepackage{bbm}
\usepackage{thmtools}
\usepackage{mathtools}
\usepackage{verbatim}
\usepackage{cancel} 
\usepackage{tikz}
\usetikzlibrary{cd}
\usetikzlibrary{calc}
\usetikzlibrary{math}
\usetikzlibrary{arrows}
\usetikzlibrary{fit}
\usetikzlibrary{positioning}
\usetikzlibrary{decorations.pathmorphing, arrows.meta}
\usetikzlibrary{shapes} 
\usetikzlibrary{snakes} 
\usetikzlibrary{matrix} 
\tikzset{Rightarrow/.style={double equal sign distance,>={Implies},->},
triple/.style={-,preaction={draw,Rightarrow}},
quadruple/.style={preaction={draw,Rightarrow,shorten >=0pt},shorten >=1pt,-,double,double
distance=0.2pt}}

\usepackage[T1]{fontenc}
 \usepackage[utf8]{inputenc}
 \usepackage{array}
 \usepackage{stmaryrd} 
 \usepackage{blkarray} 
 \usepackage{mathtools} 
 \usepackage{longtable}
 \usepackage{multicol} 
 \usepackage{xr-hyper}
 \usepackage[bookmarksdepth=2]{hyperref}
\usepackage[capitalise]{cleveref}

\usetikzlibrary{positioning} 
\usetikzlibrary{decorations.pathreplacing}
\usetikzlibrary{arrows, decorations.markings} 

\UseAllTwocells
\usepackage[
textwidth=3cm,
textsize=small,
colorinlistoftodos]
{todonotes}

\newcolumntype{D}{>{\hfil$}p{1.2cm}<{$\hfil}}

\tikzset{%
    symbol/.style={%
        draw=none,
        every to/.append style={%
            edge node={node [sloped, allow upside down, auto=false]{$#1$}}}
    }
}

\tikzset{%
scalearrow/.style n args={3}{
  decoration={
    markings,
    mark=at position (1-#1)/2*\pgfdecoratedpathlength
      with {\coordinate (#2);},
    mark=at position (1+#1)/2*\pgfdecoratedpathlength
      with {\coordinate (#3);},
    },
  postaction=decorate,
  }
}

\theoremstyle{plain}   
\newtheorem{thm}{Theorem}[section] 
\makeatletter\let\c@thm\c@thm\makeatother
\newtheorem{cor}{Corollary}[section]
\makeatletter\let\c@cor\c@thm\makeatother
\newtheorem{lem}{Lemma}[section]
\makeatletter\let\c@lem\c@thm\makeatother
\newtheorem{prop}{Proposition}[section]
\makeatletter\let\c@prop\c@thm\makeatother

\makeatletter\let\c@claim\c@thm\makeatother

\makeatletter\let\c@conjecture\c@thm\makeatother

\makeatletter\let\c@wconjecture\c@thm\makeatother

\makeatletter\let\c@thmRestate\c@thm\makeatother

\newtheorem*{unnumberedtheorem}{Theorem}

\theoremstyle{definition}

\newtheorem{defn}{Definition}[section]
\makeatletter\let\c@defn\c@thm\makeatother

\makeatletter\let\c@const\c@thm\makeatother
\newtheorem{notn}{Notation}[section]
\makeatletter\let\c@notn\c@thm\makeatother

\makeatletter\let\c@convention\c@thm\makeatother

\makeatletter\let\c@convention\c@thm\makeatother

\theoremstyle{remark}

\newtheorem{rmk}{Remark}[section]
\makeatletter\let\c@rmk\c@thm\makeatother
\newtheorem{ex}{Example}[section]
\makeatletter\let\c@ex\c@thm\makeatother

\makeatletter\let\c@observation\c@thm\makeatother

\makeatletter\let\c@warning\c@thm\makeatother

\makeatletter\let\c@digression\c@thm\makeatother

\makeatletter\let\c@answ\c@thm\makeatother

\makeatletter\let\c@answ\c@thm\makeatother

\makeatletter\let\c@aside\c@thm\makeatother

\makeatletter
\let\c@equation\c@thm
\numberwithin{equation}{section}
\makeatother

\crefname{lem}{Lemma}{Lemmas}
\crefname{thm}{Theorem}{Theorems}
\crefname{defn}{Definition}{Definitions}
\crefname{notn}{Notation}{Notations}
\crefname{const}{Construction}{Constructions}
\crefname{prop}{Proposition}{Propositions}
\crefname{rmk}{Remark}{Remarks}
\crefname{cor}{Corollary}{Corollaries}
\crefname{equation}{Display}{Displays}
\crefname{ex}{Example}{Examples}
\crefname{thmalph}{Theorem}{Theorems}
\crefname{answ}{Answer}{Answers}
\crefname{question}{Question}{Questions}

\newcommand{\cA}{\mathcal{A}}
\newcommand{\cB}{\mathcal{B}}
\newcommand{\cC}{\mathcal{C}}
\newcommand{\cD}{\mathcal{D}}

\newcommand{\cO}{\mathcal{O}}

\newcommand{\cS}{\mathcal{S}}

\newcommand{\cat}{\cC\!\mathit{at}}

\newcommand{\adch}{\mathit{ad}\mathcal C\mathit{h}}

 \newcommand{\omegacat}{\omega\cat}

 \newcommand{\comp}{\ast}

\newcommand{\inttrunc}[1]{\tau_{\leq n}^{\text{i}}}

\DeclareMathOperator{\id}{id}

\newcommand{\atom}[1]{\langle#1\rangle}

\newcommand{\tabld}[2]{\begin{pmatrix}#1^-_0 &\dots &#1^-_{#2-1}
  &#1^-_{#2}\cr\noalign{\vskip 3pt} #1^+_0 &\dots &#1^+_{#2-1}
  &#1^+_{#2}\end{pmatrix}}

\newcommand{\supp}{\mathrm{supp}}

\newcommand\pdfomega{\texorpdfstring{$\omega$}{w}}
\newcommand*\cocolon{%
        \nobreak
        \mskip6mu plus1mu
        \mathpunct{}%
        \nonscript
        \mkern-\thinmuskip
        {:}%
        \mskip2mu
        \relax
}

\AtBeginDocument{%
   \def\MR#1{}
}